\documentstyle{article}

\newcommand{\N}{\mbox{\rm I\kern-1.5pt N}}
\newcommand{\Z}{\mbox{\sf Z \hspace{-1.1em} Z}}
\newcommand{\R}{\mbox{\rm I\kern-1.5pt R}}
\newcommand{\F}{\mbox{${\rm I} \! {\rm F}$}}
\newcommand{\PP}{\mbox{${\rm I} \! {\rm P}$}}
\newcommand{\num}{${\rm n}^{\underline{\rm o}}$}
\newcommand{\dem}{\par{\em Proof\/: }\\ \noindent }
\newcommand{\findemo}{$\;\;\Box$\\}
\newfont{\frack}{eufm10}

\newtheorem{defi}{Definition}[section]
\newtheorem{lem}[defi]{Lemma}
\newtheorem{prop}[defi]{Proposition}
\newtheorem{thm}[defi]{Theorem}
\newtheorem{coro}[defi]{Corollary}
\newtheorem{nota}[defi]{Remark}

\begin{document}

\title{Rational points, genus and asymptotic behaviour in reduced 
algebraic curves over finite fields}
\author{J. I. Farr\'{a}n\thanks{J. I. Farr\'{a}n is with Dpto. Matem\'{a}tica Aplicada a la Ingenier\'{\i}a, 
E.T.S.I.I., University of Valladolid, Spain; e-mail: ignfar@eis.uva.es }}
\date{October 26, 1999}
\maketitle

\begin{abstract}

The number $A(q)$ shows the asymptotic behaviour of the quotient of 
the number of rational points over the genus of non-singular absolutely 
irreducible curves over $\F_{q}\,$. Research on bounds for $A(q)$ is 
closely connected with the so-called asymptotic main problem in Coding Theo\-ry. 
In this paper, we study some generalizations of this number for non-irreducible 
curves, their connection with $A(q)$ and its application in Co\-ding Theory. 


\end{abstract}

\section{Introduction}

For $q=p^{s}$ with $p$ prime, let ${\F}_{q}$ be the finite field with 
$q$ elements. Coding Theo\-ry is interested in the search of 
{\em asymptotically good families} of error-correcting codes, 
that is, families of codes whose relative parameters have 
a limit point in the {\em code domain} over ${\F}_{q}\,$; 
in such a family, the length of the codes cannot be upper bounded 
(see \cite{Lach} for further details). 
In the case of geometric Goppa codes, that is, codes constructed 
over non-singular absolutely irreducible algebraic curves $X$ defined 
over ${\F}_{q}\,$, Tsfasman proved in \cite{Tsf} that this problem is connected 
with the asymptotic behaviour of the quotient of the number of rational points 
over ${\F}_{q}$ by the genus of the curves, which is given over the number 
$$A(q)\doteq\limsup_{g\rightarrow\infty}\frac{N_{q}(g)}{g}$$ 
where $N_{q}(g)\doteq max\{N_{q}(X)$ : $X$ $non$-$singular$ $absolutely$ 
$irreducible$ $algebraic$ $curve$ $defined$ $over$ ${\F}_{q}$ $with$ $g(X)=g\}$ 
(see \cite{TsfVla} for further details). 
For this type of codes, the length is closely connected with the number of 
rational points, and so this number cannot be upper bounded. 
The known results about $A(q)$ can be summarized as follows: 
\begin{description}
\item[I)] $A(q)>0$ for all $q$\/. 
\item[II)] $A(q)\leq q^{\frac{1}{2}}-1$ for all $q$\/. 
\item[III)] If $q\in{\Z}^{2}$ then $A(q)=q^{\frac{1}{2}}-1$. 
\end{description}

Statement I was proved by Serre in \cite{Serre}, statement II by Drinfeld and Vladut 
in \cite{DrinVla}, and statement III independently by different people, for example 
by Garc\'{\i}a and Stichtenoth in \cite{GStich}. 

Thus a general upper bound for $A(q)$ is known, but nobody knows 
neither its exact value for all $q$ nor a general lower bound for it. 
In this paper, we give an alternative way to compute $A(q)$ by using 
non-irreducible curves under certain hypothesis, which may be useful for 
bounding $A(q)$ and so for the main asymptotic problem in Coding Theory. 
The reason of this work comes from the possibility of finding families 
of curves given by singular plane models for which it is tipically difficult 
to decide whether they are irreducible or not.

\section{$A(q)$ for non-irreducible curves}

For each $1\leq r\in{\N}$ and $q=p^{s}$ with $p$ prime, 
let ${\cal A}_{r}(q)$ be the set of 
all posible curves defined over the finite field ${\F}_{q}$ 
which are non-singular, reduced and with exactly $r$ disjoint irreducible 
components over $\overline{{\F}_{q}}$ which are all of them defined 
over ${\F}_{q}\,$, and let ${\cal A}^{r}(q)$ be the set of 
all posible curves defined over the finite field ${\F}_{q}$ 
which are non-singular, reduced and with at most $r$ disjoint irreducible 
components over $\overline{{\F}_{q}}$ which are all of them defined 
over ${\F}_{q}\,$. 

If one writes $X\in{\cal A}_{r}(q)$ as a disjoint union 
$X=X_{1}\cup\ldots\cup X_{r}\,$, 
with $X_{i}\in{\cal A}_{1}(q)$\/, one has 
$$N_{q}(X)=\sum_{i=1}^{r}N_{q}(X_{i})\;\;\;\;{\rm and}\;\;\;\;
2g(X)-2=\sum_{i=1}^{r}(2g(X_{i})-2)$$
Thus $g(X)=1-r+\sum_{i=1}^{r}g(X_{i})$ and hence 
$g(X_{i})\leq g(X)+r-1$ for all $i=1,\ldots,r$ 
(see the details in \cite{Hart}).

\begin{defi}

We define the following subsets of ${\R}$:
$${\mbox{\frack A}}_{r}(q)\doteq\{\limsup_{n\rightarrow\infty}\frac{N_{q}(X_{n})}
{g(X_{n})}\;\;:\;\;X_{n}\in{\cal A}_{r}(q)\;\;{\rm and}\;\;
\lim_{n\rightarrow\infty}N_{q}(X_{n})=\infty\}$$
$${\mbox{\frack A}}^{r}(q)\doteq\{\limsup_{n\rightarrow\infty}\frac{N_{q}(X_{n})}
{g(X_{n})}\;\;:\;\;X_{n}\in{\cal A}^{r}(q)\;\;{\rm and}\;\;
\lim_{n\rightarrow\infty}N_{q}(X_{n})=\infty\}$$

\end{defi}

\begin{nota} 

$\;$

\begin{description}

\item[a)] Notice that $\lim_{n\rightarrow\infty}N_{q}(X_{n})=\infty$
implies in both cases that $g(X_{n})$ is not upper bounded, by using 
the Hasse-Weil bound (see \cite{Weil}). So we can assume that 
$\lim_{n\rightarrow\infty}g(X_{n})=\infty$ 
and $g(X_{n})\neq0$ for all $n$\/. 

\item[b)] For each $\alpha\in{\mbox{\frack A}}_{r}(q),{\mbox{\frack A}}^{r}(q)$ there exists a 
sequence $\{X_{n}\}_{n=1}^{\infty}$ in ${\cal A}_{r}(q),{\cal A}^{r}(q)$ 
with $\lim_{n\rightarrow\infty}N_{q}(X_{n})=\infty$, 
$\lim_{n\rightarrow\infty}g(X_{n})=\infty$, $g(X_{n})\geq1$ and 
$$\lim_{n\rightarrow\infty}\frac{N_{q}(X_{n})}{g(X_{n})}=\alpha$$
being both $\{N_{q}(X_{n})\}_{n=1}^{\infty}$ 
and $\{g(X_{n})\}_{n=1}^{\infty}$ increasing. 

\item[c)] Also notice that $\mbox{\frack A}_{r}(q),\mbox{\frack A}^{r}(q)$ are 
bounded subsets of ${\R}$, since by using the Hasse-Weil bound one has 
$$\frac{N_{q}(X)}{g(X)}\leq r(1+q)+2r^{2}\sqrt{q}$$
if $X\in{\cal A}_{r}(q),{\cal A}^{r}(q)$ and $g(X)\geq1$, 
and they are lower bounded by $0$ by using the previous remarks. 

\end{description}

\end{nota}

\begin{defi}

Now define the following real numbers: 
$$A_{r}(q)\doteq sup\;{\mbox{\frack A}}_{r}(q)\;\;\;\;{\rm and}\;\;\;\;
A^{r}(q)\doteq sup\;{\mbox{\frack A}}^{r}(q)$$

\end{defi}

\begin{lem}

The sets ${\mbox{\frack A}}_{r}(q)$ and ${\mbox{\frack A}}^{r}(q)$ are 
closed for the usual topology in ${\R}$. 

\end{lem}

\dem

Take $\{\alpha_{n}\}_{n=1}^{\infty}\subseteq{\mbox{\frack A}}_{r}(q)$ (resp. 
${\mbox{\frack A}}^{r}(q)$\/) with $\lim_{n\rightarrow\infty}\alpha_{n}=\alpha$ 
and write $\alpha_{n}=\lim_{m\rightarrow\infty}\frac
{N_{q}(X_{m}^{n})}{g(X_{m}^{n})}$ as in remark 1.b. 
Then, there exists $X_{m_{1}}^{1}\equiv Y_{1}$ such that 
$\left|\frac{N_{q}(Y_{1})}{g(Y_{1})}-\alpha_{1}\right|<1$, and by induction, 
for all $n>1$ there exists $X_{m_{n}}^{n}\equiv Y_{n}$ with 
$\left|\frac{N_{q}(Y_{n})}{g(Y_{n})}-
\alpha_{n}\right|<\frac{1}{n}$ and $N_{q}(Y_{n})>N_{q}(Y_{n-1})$. 
In particular $\lim_{n\rightarrow\infty}N_{q}(Y_{n})=\infty$\/. 

Now if $\varepsilon>0$, there exist $n_{1}\in{\N}$ with 
$\frac{1}{n_{1}}<\frac{\varepsilon}{2}$ and $n_{2}\in{\N}$ with 
$|\alpha_{n}-\alpha|<\frac{\varepsilon}{2}$ for $n\geq n_{2}\,$. 
Set $n_{0}=max\{n_{1},n_{2}\}$\/; thus, if $n\geq n_{0}$ then 
$$\left|\frac{N_{q}(Y_{n})}{g(Y_{n})}-\alpha\right|\leq
\left|\frac{N_{q}(Y_{n})}{g(Y_{n})}-\alpha_{n}\right|+|\alpha_{n}-\alpha|<
\varepsilon$$
and so $\lim_{n\rightarrow\infty}\frac{N_{q}(Y_{n})}{g(Y_{n})}=\alpha$\/. 

\findemo

As a consequence, one actually has 
$$A_{r}(q)=max\;{\mbox{\frack A}}_{r}(q)\;\;\;\;{\rm and}\;\;\;\;
A^{r}(q)=max\;{\mbox{\frack A}}^{r}(q)$$

\begin{coro} $A^{1}(q)=A_{1}(q)=A(q)$ \end{coro}

\dem

The first equality is obvious. Now one has 
$$A(q)=\lim_{n\rightarrow\infty}\frac{N_{q}(X_{n})}{g(X_{n})}$$
for a sequence $\{X_{n}\}_{n=1}^{\infty}$ in ${\cal A}_{1}(q)$ with 
$N_{q}(X_{n})=N_{q}(g(X_{n}))$, $\{g(X_{n})\}_{n=1}^{\infty}$ increasing and 
$\lim_{n\rightarrow\infty}g(X_{n})=\infty$\/. But then 
$\lim_{n\rightarrow\infty}N_{q}(X_{n})=\infty$, since $A(q)>0$; 
thus $A(q)\in {\mbox{\frack A}}_{1}(q)$ and therefore $A(q)\leq A_{1}(q)$\/. 

Conversely, if one writes $A_{1}(q)=
\lim_{n\rightarrow\infty}\frac{N_{q}(X_{n})}{g(X_{n})}$ as in 
remark 1.b, since $N_{q}(X_{n})\leq N_{q}(g(X_{n}))$ one has 
$A_{1}(q)\leq\lim_{n\rightarrow\infty}\frac{N_{q}(g(X_{n}))}{g(X_{n})}\leq 
A(q)$\/. 

\findemo

\section{Properties of $A_{r}(q)$ and $A^{r}(q)$}

The following result shows some elementary properties of the numbers 
$A_{r}(q)$ and $A^{r}(q)$\/, and the proof is left to the reader.

\begin{lem}

For all $r\geq1$ one has: 
$$\begin{array}{ll}{\bf i)}\;A_{r}(q)\leq A^{r}(q) & {\bf iii)}\;A^{r}(q)>0 \\
{\bf ii)}\;A^{r}(q)\leq A^{r+1}(q) & {\bf iv}\;A_{r}(q),A^{r}(q)<\infty \end{array}$$

\end{lem}

The following formula is not so elementary, and it shows in particular 
that $A_{r}(q)>0$ for all $r\geq1$.

\begin{prop}

$A_{r}(q)=A^{r}(q)$ for all $r\geq1$. 

\end{prop}

\dem

We only have to prove that $A^{r}(q)\leq A_{r}(q)$. 
Let $X=X_{1}\cup\ldots\cup X_{s}$ in ${\cal A}_{s}(q)$ with $1\leq s<r$ and 
$g(X)\geq1$, and let $Z\in{\cal A}_{1}(q)$ with $g(Z)=1$. Then define 
$$Y=X_{1}\cup\ldots\cup X_{r}$$ 
where $X_{i}\cong Z$ for $i=s+1,\ldots,r$ and such that $X_{i}$ are 
disjoint for all $i=1,\ldots,r$. So we have got $Y\in{\cal A}_{r}(q)$ such that 
$N_{q}(X)\leq N_{q}(Y)$ and $g(X)=g(Y)\geq1$, and hence 
$$\frac{N_{q}(X)}{g(X)}\leq\frac{N_{q}(Y)}{g(Y)}$$
Now we get $A^{r}(q)\in{\mbox{\frack A}}^{r}(q)$ and write 
$A^{r}(q)=\lim_{n\rightarrow\infty}\frac{N_{q}(X_{n})}{g(X_{n})}$ 
as in remark 1.b. As in the last construction, we get a sequence 
$\{Y_{n}\}_{n=1}^{\infty}$ in ${\cal A}_{r}(q)$ such that 
$$N_{q}(X_{n})\leq N_{q}(Y_{n})\;\;\;\;{\rm and}\;\;\;\;g(X_{n})=g(Y_{n})\geq1$$

\newpage

\noindent
Thus, one has $\lim_{n\rightarrow\infty}N_{q}(Y_{n})=\infty$,  
$\lim_{n\rightarrow\infty}g(Y_{n})=\infty$ and 
$$\frac{N_{q}(X_{n})}{g(X_{n})}\leq\frac{N_{q}(Y_{n})}{g(Y_{n})}$$
for all $n\in{\N}$ and hence 
$$A^{r}(q)=\lim_{n\rightarrow\infty}\frac{N_{q}(X_{n})}{g(X_{n})}\leq
\limsup_{n\rightarrow\infty}\frac{N_{q}(Y_{n})}{g(Y_{n})}\in
{\mbox{\frack A}}_{r}(q)$$
what proofs the inequality. 

\findemo

Now in order to prove the main result of this paper 
we need some some elementary facts about limits.

\begin{lem}

Let $\{a_{n}\}_{n=1}^{\infty}$, $\{b_{n}\}_{n=1}^{\infty}$ and 
$\{c_{n}\}_{n=1}^{\infty}$ be three sequences in ${\R}$. Then:

\begin{description}

\item[i)] If $\lim_{n\rightarrow\infty}a_{n}={\pm}\infty$ then 
$a_{n}\sim a_{n}+M$, for all $M\in{\R}$.

\item[ii)] If $a_{n}\sim b_{n}$ then $\limsup_{n\rightarrow\infty}
(a_{n}c_{n})=\limsup_{n\rightarrow\infty}(b_{n}c_{n})$.

\item[iii)] If $\lim_{n\rightarrow\infty}a_{n}=
\lim_{n\rightarrow\infty}b_{n}=+\infty$ then 
$\displaystyle\limsup_{n\rightarrow\infty}\displaystyle\frac{a_{n}+M}{b_{n}+N}=
\displaystyle\limsup_{n\rightarrow\infty}\displaystyle\frac{a_{n}}{b_{n}}$. 

\item[iv)] If $a_{n}\leq b_{n}$ for $n>>0$ then 
$\limsup_{n\rightarrow\infty}a_{n}\leq\limsup_{n\rightarrow\infty}b_{n}$.  

\item[v)] If $\limsup_{n\rightarrow\infty}a_{n}<\infty$ and 
$\limsup_{n\rightarrow\infty}b_{n}<\infty$ then 
$$\limsup_{n\rightarrow\infty}(a_{n}+b_{n})\leq
\limsup_{n\rightarrow\infty}a_{n}+\limsup_{n\rightarrow\infty}b_{n}.$$ 

\end{description}

\end{lem}

\begin{thm}

$A_{r}(q)=A(q)$, for all $r\geq1$.

\end{thm}

\dem

In order to prove $A_{r}(q)\leq A_{1}(q)$ for $r\geq1$, 
we proceed by induction, and assume that $A_{r}(q)\leq A_{1}(q)$\/. 
Write $A_{r+1}(q)=\lim_{n\rightarrow\infty}\frac{N_{q}(X_{n})}{g(X_{n})}$ 
as in remark 1.b and $X_{n}=X_{n}^{1}\cup\ldots\cup X_{n}^{r+1}$ 
as a disjoint union with $X_{n}^{i}\in{\cal A}_{1}(q)$, and denote 
$$N_{q}(X_{n})=N_{n},\;\;N_{q}(X_{n}^{i})=N_{n}^{i},\;\;
g(X_{n})=g_{n},\;\;g(X_{n}^{i})=g_{n}^{i}$$

\begin{description}

\item[\underline{Case (1):}] the sequence $\{N_{n}^{i}\}_{n=1}^{\infty}$ 
is upper bounded for some $i$\/. Assumed $i=r+1$ then one has 
$N_{n}^{r+1}\leq N$ for all $n$\/. But since 
$\lim_{n\rightarrow\infty}N_{n}=\infty$, one has 
$\lim_{n\rightarrow\infty}(N_{n}^{1}+\ldots+N_{n}^{r})=\infty$ and 
$$A_{r+1}(q)=\lim_{n\rightarrow\infty}\frac{N_{n}^{1}+\ldots+N_{n}^{r}+
N_{n}^{r+1}}{g_{n}^{1}+\ldots+g_{n}^{r}+g_{n}^{r+1}-r}\leq
\limsup_{n\rightarrow\infty}\frac{N_{n}^{1}+\ldots+N_{n}^{r}+N}
{g_{n}^{1}+\ldots+g_{n}^{r}-r}=$$
$$=\limsup_{n\rightarrow\infty}
\frac{N_{n}^{1}+\ldots+N_{n}^{r}}{g_{n}^{1}+\ldots+g_{n}^{r}+1-r}\in
{\mbox{\frack A}}_{r}(q)$$
Thus $A_{r+1}(q)\leq A_{r}(q)\leq A_{1}(q)$ by the induction hypothesis. 
  
\item[\underline{Case (2):}] none of the sequences $\{N_{n}^{i}\}_{n=1}^{\infty}$ 
are upper bounded. In particular, the sequence $\{N_{n}^{r+1}\}_{n=1}^{\infty}$ 
has a subsequence with $\lim_{k\rightarrow\infty}N_{n_{k}}^{r+1}=\infty$\/; 
if $\{\sum_{i=1}^{r}N_{n_{k}}^{i}\}_{k=1}^{\infty}$ is upper bounded then the 
sequence $\{X_{n_{k}}\}_{k=1}^{\infty}$ is in Case (1) and we have finished. 
Otherwise we can get a subsequence 
$\{X_{n_{k_{l}}}\}_{l=1}^{\infty}$ such that $\lim_{l\rightarrow\infty}
(\sum_{i=1}^{r}N_{n_{k_{l}}}^{i})=\infty$ and we can 
assume that $\lim_{n\rightarrow\infty}N_{n}^{r+1}=
\lim_{n\rightarrow\infty}(\sum_{i=1}^{r}N_{n}^{i})=\infty$, 
arguing with the sequence $\{X_{n_{k_{l}}}\}_{l=1}^{\infty}$ instead of 
$\{X_{n}\}_{n=1}^{\infty}$. Then write 
$$A_{r+1}(q)=\lim_{n\rightarrow\infty}\frac{N_{n}^{1}+\ldots+N_{n}^{r}+
N_{n}^{r+1}}{g_{n}^{1}+\ldots+g_{n}^{r}+g_{n}^{r+1}-r}=
\lim_{n\rightarrow\infty}
\frac{N_{n}^{1}+\ldots+N_{n}^{r}+
N_{n}^{r+1}}{g_{n}^{1}+\ldots+g_{n}^{r}+g_{n}^{r+1}}=$$
$$=\lim_{n\rightarrow\infty}\frac{\frac{\sum_{i=1}^{r}N_{n}^{i}}
{\sum_{i=1}^{r}g_{n}^{i}}+\frac{N_{n}^{r+1}}{g_{n}^{r+1}}
\frac{g_{n}^{r+1}}{\sum_{i=1}^{r}g_{n}^{i}}}{1+\frac{g_{n}^{r+1}}
{\sum_{i=1}^{r}g_{n}^{i}}}$$  

Now one has 
$$\limsup_{n\rightarrow\infty}\frac{N_{n}^{r+1}}{g_{n}^{r+1}}\leq 
A_{1}(q)$$ and $$\limsup_{n\rightarrow\infty}\frac{\sum_{i=1}^{r}N_{n}^{i}}
{\sum_{i=1}^{r}g_{n}^{i}}=\limsup_{n\rightarrow\infty}\frac{\sum_{i=1}^{r}
N_{n}^{i}}{1-r+\sum_{i=1}^{r}g_{n}^{i}}\leq A_{r}(q)\leq A_{1}(q)$$ 
by induction hypothesis. Thus, for every $\varepsilon>0$ one has 
$$\frac{N_{n}^{r+1}}{g_{n}^{r+1}},\frac{\sum_{i=1}^{r}N_{n}^{i}}
{\sum_{i=1}^{r}g_{n}^{i}}\leq A_{1}(q)+\varepsilon$$
if $n>>0$, by the definition of upper limit, and hence 
$$A_{r+1}(q)\leq limsup_{n\rightarrow\infty}\frac{[A_{1}(q)+\varepsilon]
\left(1+\frac{g_{n}^{r+1}}{\sum_{i=1}^{r}g_{n}^{i}}\right)}
{\left(1+\frac{g_{n}^{r+1}}{\sum_{i=1}^{r}g_{n}^{i}}\right)}=
A_{1}(q)+\varepsilon$$
But since the inequality holds for all $\varepsilon>0$, we conclude 
$A_{r+1}(q)\leq A_{1}(q)$\/, what proves the theorem. 

\end{description}

\findemo

\section{Conclusion}

We may see now what happens if the number of irreducible components of 
the curves we use is not upper bound, that is, we look at the sets 
$${\cal A}^{\infty}(q)\doteq\cup_{r=1}^{\infty}{\cal A}_{r}(q)$$
$${\mbox{\frack A}}^{\infty}(q)\doteq\{\limsup_{n\rightarrow\infty}\frac
{N_{q}(X_{n})}{g(X_{n})}\;:\;X_{n}\in{\cal A}^{\infty}(q),\;
g(X_{n})\neq0\;{\rm and}\;\lim_{n\rightarrow\infty}N_{q}(X_{n})=\infty\}$$ 
and the $number$ $A^{\infty}(q)\doteq sup\;{\mbox{\frack A}}^{\infty}(q)\in\overline{\R}$. 

Firstly, we will see that $\lim_{n\rightarrow\infty}g(X_{n})=\infty$ is not 
necessary in ${\mbox{\frack A}}^{\infty}(q)$. Let $X\in{\cal A}_{1}(q)$ with 
$g(X)=1$ and $N_{q}(X)\geq1$ (note that such a curve exists for all $q\/$); 
now, for all $n\geq1$ we get $Z_{n}=X_{1}^{n}\cup\ldots\cup X_{n}^{n}$ 
disjoint union, with $X_{i}^{n}\cong X$ for all $i$ and $n$. One has 
$$N_{q}(Z_{n})\geq n\;\;{\rm and}\;\;g(Z_{n})=1\;\;{\rm for}\;\;{\rm all}\;\;n$$ 
Thus, $\lim_{n\rightarrow\infty}N_{q}(Z_{n})=\infty$ although 
$\{g(Z_{n})\}_{n=1}^{\infty}$ is constant. Furthermore, one has 
$$\limsup_{n\rightarrow\infty}\frac{N_{q}(Z_{n})}{g(Z_{n})}\geq
\limsup_{n\rightarrow\infty}n=+\infty$$ and hence $A^{\infty}(q)=
+\infty$. 

We may include the hypothesis $\lim_{n\rightarrow\infty}g(X_{n})=\infty$ in 
the definition, that is, define 
$${\mbox{\frack A}}_{\infty}(q)\doteq\{\limsup_{n\rightarrow\infty}\frac
{N_{q}(X_{n})}{g(X_{n})}\;:\;X_{n}\in{\cal A}^{\infty}(q),\;
\lim_{n\rightarrow\infty}N_{q}(X_{n})=\lim_{n\rightarrow\infty}g(X_{n})=\infty\}$$ 
and $A_{\infty}(q)\doteq sup\;{\mbox{\frack A}}_{\infty}(q)\in\overline{\R}$. 

We will see that $A_{\infty}(q)$ is not finite either; for it, take 
$Y\in{\cal A}_{1}(q)$ with $g(Y)=2$ and $N_{q}(Y)\geq1$ (also exists 
such a curve for all $q\/$) and, for $n\geq1$, define 
$Z_{n}=X_{1}^{n}\cup\ldots\cup X_{n}^{n}$ disjoint union, with 
$X_{i}^{n}\cong X$ or $X_{i}^{n}\cong Y$ for all $i$ and $n$ ($X$ as in 
example before). Note that $1\leq g(Z_{n})\leq 2n+1-n=n+1$ for all $n$. 

Now, since $1\leq\lfloor\log n\rfloor\leq n+1$ for $n>>0$, we can choose 
$X_{i}^{n}$ in $Z_{n}$ such that $g(Z_{n})=\lfloor\log n\rfloor$ for $n>>0$, 
getting $X_{i}^{n}\cong Y$ for $1\leq i\leq\lfloor\log n\rfloor-1$ and 
$X_{i}^{n}\cong X$ for $\lfloor\log n\rfloor\leq i\leq n$. Then 
$$\limsup_{n\rightarrow\infty}\frac{N_{q}(Z_{n})}{g(Z_{n})}\geq
\limsup_{n\rightarrow\infty}\frac{n}{\log n}=+\infty$$ 
and one also obtains $A_{\infty}(q)=+\infty$. 

This shows that, for calculating $A(q)$ or giving lower bounds for it, 
we must take families of curves 
with a number of irreducible components upper bounded; this hypothesis and 
the property of the irrecucible components being defined over ${\F}_{q}$ 
can be guaranteed by conditions on the branches at infinity, as we see in the following 

\newpage

\begin{coro}

Let $\{Y_{n}\}_{n=1}^{\infty}$ be a family of curves in 
${\PP}^{t_{n}}\equiv{\PP}^{t_{n}}({\F}_{q})$ defined 
over ${\F}_{q}$ and let $r$ be and integer; assume that for every $n$ 
there exists an hyperplane $H_{n}\subseteq{\PP}^{t_{n}}$ such that all of 
branches at the points in $H_{n}\cap Y_{n}$ are defined over ${\F}_{q}$ 
and they are $r$ at most. If $X_{n}$ is the normalization of $Y_{n}$ for 
every $n$ and $\lim_{n\rightarrow\infty}N_{q}(X_{n})=\infty$, then one has 
$$\limsup_{n\rightarrow\infty}\frac{N_{q}(X_{n})}{g(X_{n})}\leq A(q)$$

\end{coro}

\dem 

Every irreducible component of $Y_{n}$ over ${\F}_{q}$ has 
at least one point in $H_{n}$; by hypothesis, the branches at such a point 
are defined over ${\F}_{q}$, then the irreducible components of $Y_{n}$ 
(and so the ones of $X_{n}$) are defined over ${\F}_{q}$ and the are $r$ 
at most. Hence the statement follows from the above theorem. 

\findemo

\end{document}